\documentclass{elsart}

\usepackage{amsfonts}
\usepackage{amstext}
\usepackage{amssymb,amsmath}

\def\R{\mathbb R}
\def\EE{\mathbb E}\def\PP{\mathbb P}
\def\0{\mathbf 0}
\def\1{\mathbf 1}
\def\a{\mathbf a}

\def\x{\mathbf x}

\def\z{\mathbf z}
\def\y{\mathbf y}
\def\h{\mathbf h}
\def\X{\mathbf X}

\def\A{\mathbf A}

\def\J{\mathbf J}
\def\I{\mathbf I}
\def\EEE{\mathbf E}
\def\PPP{\mathbf P}

\def\W{\mathbf W}
\def\BB{\mathbf B}

\def\eps{\varepsilon}

\def\rk{\text{rank}\,}

\def\dim{\text{dim}}

\begin{document}

\begin{frontmatter}
\title {Testability of minimum balanced multiway cut densities}
\author[bmemat]{Marianna Bolla\corauthref{ca}}, 
\ead{marib@math.bme.hu}
\corauth[ca]{Corresponding author. Research supported in part by the Hungarian 
National Research Grants OTKA 76481 and OTKA-NKTH 77778.}
\author[bmemat]{Tam\'as K\'oi},
\ead{koitomi@math.bme.hu}
\author[szeged]{Andr\'as Kr\'amli}
\ead{kramli@informatika.ilab.sztaki.hu}

\address[bmemat]{Institute of Mathematics,
Budapest University of Technology and Economics}
\address[szeged]{Bolyai Institute, University of Szeged}

\begin{abstract} 
Testable weighted graph parameters and equivalent notions of testability 
are investigated based on~\cite{LovI}. 
We prove that certain balanced minimum multiway cut densities are testable.
Using this fact, quadratic programming techniques are applied
to approximate some of these quantities.
 The problem is related to  cluster analysis  and
statistical physics. 
Convergence of special noisy graph 
sequences is also discussed.
\end{abstract}


\end{frontmatter}

\section{Introduction}

A typical problem of contemporary cluster analysis  is to find 
relatively small number of homogeneous groups of  data 
that do not differ significantly in size.
To make inferences on the separation that can be achieved for a given number
of clusters, some types  of minimum cut densities are investigated.

In a fairly general setup of~\cite{LovI}, 
the objects to be classified are vertices of
a weighted graph whose edges and vertices both have nonnegative, real weights.
Edge-weights are similarities between the vertices normalized in such a way 
that 0 is the minimum 
and 1 is the maximum similarity,
 while vertex-weights reflect individual values of the vertices.
Classical (simple) graphs have vertex-weights all equal to 1 and 
edge-weights 0 or 1.

For given number $n$ of vertices and for a fixed integer $0<q\le n$ we define 
three types of minimum $q$-way cut densities, each being the minimum of the 
weight-sum of between-cluster edges, occasionally  adjusted with a factor 
characterizing  within-cluster densities, over all or over balanced 
$q$-partitions of the vertices. 
The limit of these densities is considered as $n\to \infty$. 
If this limit exists for any convergent graph sequence, 
we say, that the $q$-way cut density in question is a testable graph parameter.
In fact, the subsequent terms of such a convergent graph sequence $(G_n )$
 become more and more 
 similar in their global structure, which fact can be formulated in terms of
convergence of the homomorphism densities of injective maps $F\to G_n$ for any 
simple graph $F$.

Hence, testable parameters measure statistical properties of a large graph
that are indifferent to minor changes in the edge- and vertex-weights. 
It will be proved that certain balanced $q$-way cut densities are testable.
To this end, notions of testability are extended to weighted graphs, and we 
prove equivalent statements of testability by means of large deviation results
of Lov\'asz and coauthors~\cite{LovI}. 
Roughly speaking, these propositions state that if a
smaller simple graph is selected -- by an appropriate randomization --
 based on a large weighted graph, 
the testable parameter of the randomized one is very close to that of the 
whole graph with high probability.

The organization of the paper is as follows.
In Section~\ref{pre}, 
notion of a convergent graph sequence and that of a graphon
is introduced based on~\cite{LovI}. 
In Section~\ref{koi}, equivalent statements of testability
are discussed for weighted graphs.
In Section~\ref{cut}, testability of different kinds of minimum multiway
cut densities is investigated. For non testable ones counterexamples are 
presented, while for testable ones theorems of~\cite{LovII} based on 
statistical physics are applied.
In Section~\ref{prob},  continuous extensions of testable weighted
graph parameters to graphons  are constructed that gives rise to a quadratic 
programming task.
In Section~\ref{noise}, 
special graph sequences (blown up structures burdened with a
very general kind of noise) are analyzed utilizing the fact that the cut-norm
of a so-called Wigner-noise tends to zero as its size tends to infinity.

\section{Preliminaries}\label{pre}


Let $G=G_n$ be a weighted graph on the vertex set  
$V(G)=\{ 1,\dots ,n \} =[n]$ and edge set $E(G)$. Both the edges and vertices
 have weights: the edge-weights are  pairwise similarities 
$\beta_{ij} =\beta_{ji} \in [0,1]$,
$i,j \in [n]$, 
while the vertex-weights
$\alpha_{i} >0$ $(i\in [n])$ indicate relative significance of the 
vertices. It is important that the edge-weights are nonnegative (zero means
no connection at all), the normalization into the [0,1] interval is
for the sake of treating them later as  probabilities for random sampling.
Let ${\cal G}$ denote the set of all such weighted graphs.

The {\it volume} of $G \in {\cal G}$ is defined by
$\alpha_G =\sum_{i=1}^n \alpha_i$, while that of the vertex-subset $T$ by
$\alpha_{T} =\sum_{i\in T} \alpha_i$.
Further,
$$
 e_G (S,T) =\sum_{s\in S} \sum_{t\in T} \alpha_s \alpha_t \beta_{st} 
$$
denotes the {\it weighted cut} between the (not necessarily disjoint) 
vertex-subsets $S$ and $T$.

Lov\'asz and coauthors~\cite{LovI}  
define the homomorphism density between the simple graph 
$F$ (on vertex set $V(F) =[k]$) and  the above weighted graph $G$.
With the notations
$$
 \alpha_{\Phi } =  \prod_{i=1}^k \alpha_{\Phi (i)} , \qquad
 {\text{inj}}_{\Phi} (F,G) =\prod_{ij\in E(F)} \beta_{\Phi (i) \Phi (j)} 
$$
the homomorphism density between $F$ and $G$ is defined by
\begin{equation}\label{hom}
 t(F,G) = \frac1{(\alpha_G )^k} \sum_{\Phi : V(F)\to V(G)} 
  \alpha_{\Phi } \cdot  {\text {inj}}_{\Phi} (F,G) .
\end{equation}
For a simple graph $G$, $t(F,G)$ is the probability that a random map
$V(F) \to V(G)$ is a homomorphism. Similarly, $t_{inj} (F,G)$ and 
 $t_{ind} (F,G)$ are defined  in such a way that for a simple $G$,
they are the probabilities that a random injective map
$V(F) \to V(G)$ is adjacency preserving and results in an induced subgraph
of $F$ in $G$, respectively. With the notation
$$
 {\text{ind}}_{\Phi} (F,G) =\prod_{ij\in E(F)} \beta_{\Phi (i) \Phi (j)}
 \prod_{ij\in E(\bar F)} (1-\beta_{\Phi (i) \Phi (j)} )  ,
$$
let
$$
 t_{\text{inj}}(F,G) := \frac1{k! (\alpha )_k} \sum_{\Phi {\text{ inj.}}} 
  \alpha_{\Phi }  \cdot {\text{inj}}_{\Phi} (F,G)
$$
and
$$
 t_{\text{ind}}(F,G) := \frac1{k! (\alpha )_k} \sum_{\Phi {\text{ inj.}}} 
  \alpha_{\Phi } \cdot  {\text{ind}}_{\Phi} (F,G) ,
$$
where $(\alpha)_k$ denotes the $k$th elementary symmetric polynomial
of $\alpha_1 ,\dots ,\alpha_n$.
Latter one resembles to the likelihood function of taking a sample 
-- that is a simple graph on
$k$ vertices -- from the weighted graph $G$ in the following way: 
 $k$ vertices are chosen {\it with replacement}
with respective probabilities $\alpha_i /\alpha_G$ $(i=1,\dots ,n)$.
Given the vertex-subset $\{ \Phi (1) ,\dots ,\Phi (k) \}$,
the edges come into existence conditionally independently, with
probabilities of the edge-weights. 
Such a random graph is denoted by $\xi (k,G)$. Obviously,
$$
 \PP (\xi (k,G)=F) = \frac1{(\alpha_G )^k} \sum_{\Phi :V(F)\to V(G)} 
  \alpha_{\Phi }   {\text{ind}}_{\Phi} (F,G) ,
$$
since we may get back $F$,  even if $\Phi$ is not injective.
As most  maps into a large graph are injective, the above probability is
very close to $t_{\text{ind}}(F,G)$, and  $t(F,G)$ is very close to 
$t_{\text{inj}}(F,G)$. Further,   $t_{\text{ind}}(F,G)$
has a well-defined relation to $t_{inj}(F,G)$ that will be formulated in
Section 3. In the sequel only the $k \ll n$ case 
 makes sense, and this is the situation we need: $k$ is kept
fixed, while $n$ tends to infinity.

\begin{defn}\label{conv}
We say that the weighted graph sequence $(G_n )$ is (left-)convergent, 
if the sequence $t(F, G_n )$ converges for any simple graph $F$ ($n\to\infty$).
\end{defn}
As other kinds of convergence are not discussed here, in the sequel the
word left will be omitted, and we simply use convergence.

Authors in~\cite{LovI} also construct the limit object that is a symmetric, 
bounded, 
measurable function $W:[0,1]\times [0,1]\to \R$ and they call it {\it graphon}.
Let $\cal W$ denote the set of these functions.
The interval [0,1] corresponds to the vertices and the values 
$W(x,y) =W(y,x)$ to the edge-weights. 
In view of the conditions imposed on the edge-weights, the range
is also the [0,1] interval. The set of symmetric, measurable  functions 
$W:[0,1]\times [0,1]\to [0,1]$ is denoted by ${\cal W}_{[0,1]}$.
The  stepfunction graphon  $W_{G}\in {\cal W}_{[0,1]}$ 
is assigned to the weighted graph $G\in {\cal G}$ in the following way: 
the sides of the unit square are divided into intervals
$I_1 , \dots ,I_n$ of lengths $\alpha_1 /\alpha_G ,\dots ,\alpha_n /\alpha_G$,
and over the rectangle $I_i \times I_j$ the stepfunction takes on the
value $\beta_{ij}$. 

The so-called {\it cut-distance} between the graphons $W$ and $U$ is
\begin{equation}\label{tav}
 \delta_{\square}  (W,U) =
\inf_{\nu}  \| W -U^{\nu} \|_{\square} 
\end{equation}
where the {\it cut-norm} of the graphon $W$ is defined by
$$
 \| W  \|_{\square} =\sup_{S,T \subset [0,1]} 
 \left | \iint_{S\times T}  W(x,y) \, dx \, dy \right | ,
$$
and the infimum in~(\ref{tav}) is taken over all  measure preserving bijections
$\nu : [0,1] \to [0,1]$, while $U^{\nu}$ denotes the transformed $U$ after 
performing the same measure preserving bijection $\nu$ on both sides of the
unit square. 
An equivalence relation is defined over the set of graphons: two graphons
belong to the same class if they can be transformed into each other by
a measure preserving map, 
i.e., their  $\delta_{\square}$-distance is zero.
In the sequel, we consider graphons modulo measure preserving maps,
and under graphon we understand the whole
equivalence class. By Theorem 5.1 of~\cite{LovSzeg2},
the classes of ${\cal W}_{[0,1]}$ form a compact metric space with the
$\delta_{\square} $ metric.

We will intensively use the following reversible relation between 
convergent weighted graph sequences and
graphons.
\begin{thm}\label{lim}
(Corollary 3.9 of~\cite{LovI}). For any convergent sequence $(G_n )$ of
weighted graphs with uniformly bounded edge-weights there exists a graphon
such that $\delta_{\square} (W_{G_n} , W) \to 0$. Conversely, any graphon
$W$ can be obtained as the limit of a sequence of
weighted graphs with uniformly bounded edge-weights. The limit of a convergent 
graph sequence is essentially unique: If $G_n \to W$, then also $G_n \to W'$
for precisely those graphons $W'$ for which $\delta_{\square} (W,W')=0$.
\end{thm}

Authors of~\cite{LovI} also define the $\delta_{\square}$-distance of two
weighted graphs and that of a graphon and a graph.
Without going into details, we just cite the following facts:
for the weighted graphs $G$, $G'$, and for the graphon $W$
$$
\delta_{\square} (G,G') =\delta_{\square} (W_G ,W_{G'})  \quad \text{and}
\quad \delta_{\square} (W,G) =\delta_{\square} (W,W_{G}) .
$$  
They  prove (Theorem 2.6) that a sequence of
weighted graphs with uniformly bounded edge-weights is convergent
if and only if it is a Cauchy sequence in the metric $\delta_{\square}$.

A simple graph on $k$ vertices can be sampled based on $W$ in the
following way: $k$ uniform random numbers, $X_1 ,\dots ,X_k$ are generated 
on [0,1] independently. Then we connect the vertices corresponding 
to $X_i$ and $X_j$ with probability $W(X_i ,X_j )$.
For the so obtained simple graph $\xi (k, W)$  
the following large deviation result is proved.
\begin{thm}\label{deviation} (Theorem 4.7 of~\cite{LovI}, part (ii)).
Let $k$ be a positive integer and $W\in {\cal W}_{[0,1]}$ be a graphon.
Then with probability at least $1-e^{-k^2/(2\log_2 k )}$, we have
\begin{equation}
 \delta_{\square}  (W, {\xi (k,W)} ) \le \frac{10}{\sqrt {\log_2 k}} .
\label{dev}
\end{equation}
\end{thm}
Fixing $k$, the inequality~(\ref{dev}) holds uniformly for any 
graphon $W \in {\cal W}_{[0,1]}$, especially for $W_G$. 
Further, the sampling from $W_G$
is identical to the previously defined sampling with replacement from $G$,
that is  $\xi (k,G) =\xi (k,W_G )$. In fact, this argument is relevant
in the $k\le |V(G)|$ case.

\section{Testable weighted graph parameters}\label{koi}

A function $f:G\to \R$ is called a {\it graph parameter}
 if it is invariant under
isomorphism. In fact, a graph parameter is a statistic evaluated on the graph,
 and hence,  we are interested in weighted 
graph parameters that are not sensitive to minor changes in the weights of
the graph. 

The testability  results of~\cite{LovI} for 
simple graphs remain valid if we consider weighted
graph sequences $(G_n )$  with {\it no dominant vertex-weights}, that is 
$$
  \max_i \frac{\alpha_i (G_n)}{\alpha_{G_n}} \to 0 , \qquad n\to\infty .
$$

\begin{defn}\label{test}
A weighted graph parameter $f$ is testable if for every $\eps >0$
there is a positive integer $k$ such that if $G\in {\cal G}$ satisfies
$$\max_i \frac{\alpha_i (G)}{\alpha_G} \le \frac1{k},
$$
then 
\begin{equation}\label{testability}
 \PP (| f(G) - f( \xi (k,G) )| > \eps ) \le \eps ,
\end{equation}
where $\xi (k,G)$ is a random simple graph on  $k$ vertices
selected randomly from $G$ with replacement as described in Section 2.
\end{defn}

Consequently, such a graph parameter can be consistently estimated 
based on a fairly large sample. As the randomization depends only on the
$\alpha_i (G) / \alpha_{G}$ ratios, it is not able to distinguish
between weighted graphs whose vertex-weights differ only in a constant
factor. Thus, a testable weighted graph parameter is invariant under scaling 
the vertex-weights.
Now, we introduce some equivalent statements of the testability, indicating
that a testable parameter depends continuously on the whole graph.
This is the generalization of Theorem 6.1 of~\cite{LovI}
applicable for simple graphs.

\begin{thm} \label{ekvi}
For the weighted graph parameter $f$ the following are equivalent:
\begin{itemize}
\item[{(a)}] $f$ is testable.
\item[{(b)}]
For every $\eps >0$ there is a positive integer $k$ such that for every
weighted graph $G\in {\cal G}$ satisfying the node-condition
$\max_i {\alpha_i (G)}/{\alpha_G} \le 1/{k}$,
$$
 | f(G) - \EE (f( \xi (k,G) ) ) | \le \eps  .
$$
\item[{(c)}]
For every convergent weighted graph sequence $(G_n)$ with
$\max_i {\alpha_i (G_n )}/{\alpha_{G_n} } \to 0$,
$f(G_n )$ is also convergent ($n\to\infty$).
\item[{(d)}] $f$ can be extended to graphons such that the graphon functional
$\tilde{f} $ is continuous in the cut-norm and
$\tilde{f} (W_{G_n} ) -f(G_n ) \to 0$, whenever
$
  \max_i {\alpha_i (G_n ) }/{\alpha_{G_n} } \to 0 
$ ($n\to\infty $).
\item[{(e)}] For every $\eps >0$ there is an   $\eps_0 >0$ real and an 
$n_0 >0$ integer such that if
$G_1, G_2$ are weighted graphs satisfying
$\max_i {\alpha_i (G_1)}/{\alpha_{G_1} } \le 1/{n_0}$,
$\max_i {\alpha_i (G_2)}/{\alpha_{G_2} } \le 1/{n_0}$,
and  $\delta_{\square} ({G_1} ,{G_2} )<\eps_0$, then
$|f(G_1)-f(G_2 )| <\eps$.
\end{itemize}
\end{thm}

To prove the theorem we need three lemmas that are partly generalizations of
results in~\cite{LovI} stated for simple  graphs.
\begin{lem}\label{l1}
If $(G_n )$ is a weighted graph sequence with no dominant vertex-weights, then
for any simple graph $F$
$$ 
 | t (F, G_n ) - t_{inj} (F, G_n ) | \to 0, \qquad n\to \infty .
$$
\end{lem}

\begin{pf}
The difference between  between $ t (F, G_n )$ and $t_{inj} (F, G_n )$ is
essentially 
obtained by the summation in~(\ref{hom}) over the non injective maps. 
As the sum is increased
if we take the non-zero $\beta_{ij}$'s 1, it suffices to prove that
$$
 \sum_{\Phi \text{ non inj.}} \frac{\alpha_{\Phi} (G_n)}{\alpha_{G_n}^k}\to 0 ,
$$
where the left hand side is  the probability that there are repetitions in the
vertices. 
As 
$
  \max_i {\alpha_i (G_n ) }/{\alpha_{G_n} } < \eps_n 
$,  this probability is less than
$$
 1-1\cdot (1-\eps_n ) \dots  (1- (k-1)\eps_n ) \le 1- (1-(k-1)\eps_n )^k
$$
that tends to 0, since  $\eps_n \to 0$, provided  $k$ is fixed.
\end{pf}

\begin{lem}\label{l2}
$$
 t_{\text{inj}} (F,G) =\sum_{F'\supseteq F} t_{\text{ind}} (F',G) \quad
 \text{and} \quad 
 t_{\text{ind}} (F,G) =\sum_{F'\supseteq F} 
 (-1)^{|E(F') \setminus E(F)|} t_{\text{inj}} (F',G)
$$
where $F'$ is a simple super-graph of $F$ (on the same vertex-set, and
edge-set containing the edge-set of $F$).
\end{lem}

\begin{pf}
To verify the first statement, it suffices to prove that for any injective map
 $\Phi : V(F)\to V(G)$,
\begin{equation}\label{egy}
 {\text{inj}}_{\Phi} (F,G) =\sum_{F'\supseteq F} {\text{ind}}_{\Phi} (F',G) .
\end{equation}
Suppose that $\ell$ edges are missing from $F$ to be 
a complete graph. 
As $F'\supseteq F$ is a super-graph of $F$, it can easily be seen that 
${\text{ind}}_{\Phi} (F',G)$ contains the multiplicative factor
 ${\text{inj}}_{\Phi} (F ,G)$. 
Hence,  the right hand side of~(\ref{egy}) can be factorized as 
${\text{inj}}_{\Phi} (F,G) \cdot S_{\ell}$, where $S_{\ell}$ depends on
$F$, $G$, and $\Phi$, but for the sake of simplicity we omit these ones. 
We show -- by reverse induction on the number 
of edges -- that $S_{\ell}=1$. 
If $F$ is a  complete graph on $k$ vertices, then by 
the definition of  ${\text{inj}}_{\Phi}$ and ${\text{ind}}_{\Phi}$, 
$S_{0} = 1$. 
If $j$ edges are missing from $F$ to be a complete graph,
denote by  $\beta_1 ,\dots ,\beta_j$ the weights of their
$\Phi$-images. Further, 
denote by $\sigma = (\sigma_{1}, \dots ,\sigma_{j})$ a $\{ 0,1 \}$ 
sequence of length $j$ and $a_{i}^{(0)} (x):=x$, $a_{i}^{(1)} 
(x):=1-x$, $i=1, \dots, j$. With this notation 
$S_j = \sum_{\sigma \in \{ 0,1 \} ^j} 
\prod_{i=1}^{j} a_{i}^{(\sigma_{i})} (\beta_i )$. 
Coupling the sequences which differ 
only in the first coordinate, $S_j$ reduces to $S_{j-1}$, etc.

By inclusion-exclusion, the second statement also follows.
\end{pf}

\begin{lem}\label{l3} (Lemma 5.3 of~\cite{LovI}).
Let $(G_n )$ be a sequence of
weighted graphs with uniformly bounded edge-weights, and no dominant 
vertex-weights. If $\delta_{\square} (U,W_{G_n}) \to 0$ for some
$U\in {\cal W}$, then  the graphs in the sequence $(G_n )$ can be
relabeled in such a way that the resulting sequence $(G'_n )$ of labeled 
graphs converges to $U$ in the cut-norm: $\| U- W_{G'_n } \|_{\square} \to 0$.
\end{lem}

Now, we are able to prove the main theorem (Theorem~\ref{ekvi}).
\begin{pf} 
The idea of the proof is analogous to that of Theorem 6.1 of~\cite{LovI}. 

First we prove that (a),(b),(c),(e) are equivalent:

$(a) \Rightarrow (b)$: The statement is obvious, as due to the
boundedness of $f$, (\ref{testability})
implies that the difference is small on average.

$(b) \Rightarrow (c)$: Let $(G_n)$ be a convergent sequence of weighted 
graphs with no dominant vertex-weights. Let $\varepsilon > 0$ be 
arbitrary, and $k$ is chosen corresponding to 
$\varepsilon$ as in statement (b). If 
$n$ is large enough, then $| f(G_n) - \mathbb{E} (f( \xi (k,G_n) ) ) | \le 
\varepsilon$. On the other hand, by the definition of convergence it 
follows that 
$t(F,G_n)$ is convergent for all simple graphs $F$ on $k$ vertices. Using 
Lemmas~\ref{l1}, \ref{l2}, $t_{\text{ind}} (F,G_n)$ tends to a limit 
value denoted by $t_{\text{ind}} (F)$.  
This means that $\PP \left( \xi (k,G_n) = F 
\right) \rightarrow t_{\text{ind}} (F)$ and so $$ \EE (f( \xi (k,G_n) ) 
) \rightarrow \sum_{F: \, |V(F)| =k} t_{\text{ind}}(F) \cdot f(F) = 
a_k, $$ 
since the number of simple graphs on $k$ vertices is finite. In summary,
$$
|f(G_n) - a_k| \le |f(G_n) - \EE (f( \xi (k,G_n) ) ) | +
|\EE (f( \xi (k,G_n) ) ) - a_k| \le 2 \varepsilon
$$ 
provided $n$ is large enough. 

$(c) \Rightarrow (e)$: Suppose that (e) does not hold. In this case 
there exist $\eps > 0$, further sequences $(G_n)$ and $(G_{n}^{'})$  of 
weighted graphs, such that the dominant vertex-weights of both sequences 
tend to $0$, $\delta_{\square} (G_n ,G_{n}^{'} ) \rightarrow 0$, and 
$|f(G_n) - f(G_{n}^{'}) | \ge \varepsilon$. Using the compactness of 
${\cal W}_{[0,1]}$ we can assume that both sequences are convergent. For this 
reason, the merged sequence $G_1, G_{1}^{'}, G_2 , G'_2 ,\dots$ is also 
convergent.  For the above merged sequence, by (c),   the sequence 
 $f(G_1 ), f(G_{1}^{'}), f(G_2 ) , f(G'_2 ) ,\dots$ 
is  covergent, that  contradicts to 
$|f(G_n) - f(G_{n}^{'})| \ge \eps$.

$(e) \Rightarrow (a)$: Suppose that (a) does not hold. In this case 
there exist $\varepsilon >0$ and a sequence $(G_n)$ such that $\max_i 
\frac{\alpha_i (G_n)}{\alpha_{G_n}} \le \frac{1}{n}$, and with 
probability at least $\varepsilon$ the inequality $| f(G_n) - f( \xi 
(n,G_n) )| > \varepsilon$ holds for all $n$. To this $\varepsilon$ 
choose the corresponding $n_0$ and $\varepsilon_{0}$ 
as in the statement (e). Furtheron,
because of (\ref{dev}), the sequence $\delta_{\square} (G_n , 
\xi (n,G_n) )$ tends to $0$ in probability. In particular,  $\PP 
\left( \delta_{\square} (G_n , \xi (n,G_n) ) < \varepsilon_0 \right) \ge 
1 - \frac{\varepsilon}{2}$. Using the definition of $\varepsilon_0$ and 
$n_0$ we get that $\PP \left( | f(G_n) - f( \xi (n,G_n) )| < 
\varepsilon \right) \ge 1 - \frac{\varepsilon}{2}$. This contradicts to 
the fact that with probability at least $\varepsilon$ the opposite is 
true.

Now we prove that the statement (d) is also equivalent to the 
testability.

$(c),(e) \Rightarrow (d)$: Let $W \in {\cal W}_{[0,1]}$ be an arbitrary 
graphon. By Theorem~\ref{lim} we can find a sequence $(G_n )$ of 
weighted graphs with no dominant vertex-weights tending to $W$. Let 
${\tilde f}(W)$ be the limit of $f(G_n)$. Because of (c) the limit 
exists, and due to the statement (e) this definition is correct.
First we prove the continuity. Let $\varepsilon >0$ be arbitrary. Using 
the statement (e), to $\frac{\varepsilon}{3}$ we assign the corresponding
 $\eps^{'}$  and $n_0$. We show, that $||W - W^{'}||_{\square} \le 
\frac{\varepsilon^{'}}{3}$ implies $| \Tilde{f}(W) - \Tilde{f}(W^{'}) | 
\le \varepsilon$. For this purpose let $(G_n)$ be a sequence of 
weighted graphs with 
no dominant vertex-weights tending to $W$. We can choose a $G$ from $(G_n)$ 
such 
that the dominant vertex-weight of $G$ is smaller than $\frac{1}{n_0}$;
further, 
$\delta_{\square} (G,W) < \frac{\varepsilon^{'}}{3}$ and $|f(G) - 
\Tilde{f} (W) | \le \frac{\varepsilon}{3}$. Similarly, we can choose a $G^{'}$ 
from the sequence $(G'_n)$ tending to $W'$ with analogous properties.
In this case $\delta_{\square} (G, G^{'}) \le 
\delta_{\square} (G, W) + \delta_{\square} (W, W^{'}) + \delta_{\square} 
(W^{'}, G^{'}) \le \varepsilon^{'}.$ By  (e), $| f(G) - f(G^{'} | \le 
\frac{\varepsilon}{3}$, and hence, 
$$
| {\tilde f}(W) - {\tilde f}(W^{'}) | \le | {\tilde f}(W) - f(G) | + | f(G) 
| - f(G^{'}) | + | f(G^{'}) - {\tilde f}(W^{'}) | \le \varepsilon.
$$
It remains to show that
$\left| {\tilde f}(W_{G_n}) - f(G_n) \right| \to 0$, whenever 
$|V(G_n )| \to\infty$ with no dominant vertex-weights.
On the contrary,
suppose that  there exists a sequence $(G_n)$ with no dominant 
vertex-weights such that ${\tilde f} (W_{G_n}) - f(G_n)$ does not tend to 
$0$. For the sake of simplicity we can assume that for some 
$\varepsilon$: $\left| {\tilde f}(W_{G_n}) - f(G_n) \right| > \varepsilon$ 
for all $n$. We can also assume that $(G_n)$ converges to some graphon 
$W$ in the $\delta_{\square}$ metric. By Lemma~\ref{l3}, there is 
a sequence $(G_{n}^{'})$ isomorphic to $(G_n)$ such that $||W_{G_{n}^{'}} - 
W||_{\square} \rightarrow 0$. Using the statement (c), 
$\lim_{n\to\infty} f(G_{n})= \lim_{n\to\infty} f(G_{n}^{'}) =
{\tilde f} (W)$. In addition, ${\tilde f}$ is continuous, and
for this reason, 
$\lim_{n\to\infty } {\tilde f} (W_{G_{n}})=
\lim_{n\to\infty } {\tilde f} (W_{G'_{n}})= 
{\tilde f} (W)$. But this is a  contradiction.

$(d) \Rightarrow (c)$: Let $(G_n)$ be a convergent sequence of weighted
graphs with 
no dominant vertex-weights. Let $W$ be its limit. So $\delta_{\square} 
(W_{G_n}, W) \rightarrow 0$. In this way, by Lemma~\ref{l3}, we 
can relabel $(G_n)$ into $(G^{'}_{n})$ in such a way that $||W_{G^{'}_{n}} 
- W||_{\square} \rightarrow 0$. Therefore,  using the continuity of 
${\tilde f}$ we get ${\tilde f} (W_{G^{'}_{n}}) - {\tilde f} (W) 
\rightarrow 0$. Since $f(G_n) - {\tilde f} (W_{G^{'}_{n}}) = f(G_n) - 
{\tilde f} (W_{G_{n}})$, the last term tends to $0$ because of the 
statement (d). Thus, $f(G_n )$ is convergent.
\end{pf}

\begin{rem} The original testability theorem for simple graphs in~\cite{LovI}
 was formulated in terms of sampling without replacement. In the most 
important case, when the size of the sample is small compared
to the size of the underlaying graph, the two sampling methods are 
approximately the same. Usually, this is the case in practical
applications. In our definition of the testability of a weighted graph
parameter we use sampling with replacement, but the testability could be
defined by any randomization for which a
large deviation result similar to that of Theorem~\ref{deviation} holds. 
However,  
equivalent statements (c), (d), (e) do not depend on the randomization,
and we may expect their equivalence to statemants (a),(b) under an
appropriate sampling with likelihood function strongly  connected to
$t_{\text{ind}} (F,G)$ and satisfying (\ref{dev}).
\end{rem}

\section{Balanced multiway cuts}\label{cut}

Lov\'asz and coauthors~\cite{LovII} proved the testability of the 
maximum cut density. The minimum
cut density is somewhat different. E.g., if a single vertex is loosely
connected to a dense part, the minimum cut density of the whole graph is
small, however, randomizing a smaller sample, with high probability, it
will come from the dense part with a large  minimum cut density.
 
To prove the testability of certain balanced minimum multiway cut densities
we use the notions of statistical physics in the same way as in~\cite{LovII}.
Most of these notions are self-explanatory. However, to be self-contained,
we included some definitions for clarification together with
the notion  of a factor graph.

Let $G\in {\cal G}$ be a weighted graph on $n$ vertices 
with vertex-weights $\alpha_1 ,\dots ,
\alpha_n$ and edge-weights $\beta_{ij}$'s.
Let $q\le n$ be a fixed positive integer, and 
${\cal P}_q $ denote the set of $q$-partitions
$P=(V_1 ,\dots ,V_q )$ of the vertex set $V$.  The non-empty, disjoint 
vertex-subsets sometimes are referred to as clusters or states.
The {\it factor graph} or $q$-{\it quotient} of $G$ with respect to the 
$q$-partition $P$ is
denoted by $G/P$ and it is defined as the weighted graph on $q$ vertices
with vertex- and edge-weights
$$
 \alpha_i (G/P) =\frac{\alpha_{V_i}}{\alpha_G} \quad (i\in [q]) \quad
 \text{and} \quad  
\beta_{ij} (G/P) =\frac{e_G (V_i ,V_j )}{\alpha_{V_i} \alpha_{V_j} } 
 \quad (i,j \in [q]),
$$
respectively. Let ${\hat {\cal S}_q} (G)$ 
denote the set of all $q$-quotients of $G$.
The Hausdorff distance between  ${\hat {\cal S}_q} (G)$ and 
${\hat {\cal S}_q} (G')$ is
defined by
$$
 d^{\text{Hf}} ({\hat {\cal S}_q} (G) , {\hat {\cal S}_q} (G') ) = \max \{
 \sup_{H \in{\hat {\cal S}}_q (G)} \inf_{H' \in {\hat {\cal S}}_q (G')} 
d_1 (H,H' ) \, , \,
  \sup_{H' \in {\hat {\cal S}}_q (G')} \inf_{H \in {\hat {\cal S}}_q (G)} 
d_1 (H,H' ) \}  ,
$$
where
$$
 d_1 (H,H' )= \sum_{i,j\in [q]} \left| 
 \frac{\alpha_i (H) \alpha_j (H) \beta_{ij} (H)}{\alpha_H^2} -
 \frac{\alpha_i (H') \alpha_j (H') \beta_{ij} (H')}{\alpha_{H'}^2} \right| +
\sum_{i\in [q]} \left| 
 \frac{\alpha_i (H) }{\alpha_H} -
 \frac{\alpha_i (H')}{\alpha_{H'}} \right|
$$
is the $l_1$-distance between two weighted graphs $H$ and $H'$ on the same
number of vertices. Here especially, $H$ and $H'$ are factor graphs, 
and hence, $\alpha_H =\alpha_{H'} =1$, 
therefore the denominators can be omitted.

Given the real symmetric $q\times q$ matrix $\J$ and the vector $\h \in \R^q$,
the partitions $P\in {\cal P}_q$ also define a spin system on the weighted
graph $G$. The so-called {\it ground state energy} 
of such a spin configuration is
$$
{\hat {\cal E}_q } (G ,\J ,\h ) =-\max_{P\in {\cal P}_q } \left(
 \sum_{i\in[q]} \alpha_i (G/P) h_i + 
 \sum_{i,j\in[q]} \alpha_i (G/P) \alpha_j (G/P) \beta_{ij} (G/P) J_{ij}\right).
$$
Here $\J$ is the so-called coupling-constant matrix, where $J_{ij}$ represents 
the strength of interaction between states $i$ and 
$j$, and $\h$ is the magnetic field. They carry physical meaning. We shall
use only special $\J$ and $\h$, especially $\h =\0$.

Sometimes, we need balanced $q$-partitions to regulate the proportion of the
cluster volumes. A slight balancing 
between the cluster volumes is  achieved by fixing a positive real number
$c$ ($c\le 1/q$). Let
${\cal P}^c_q$ denote the set of $q$-partitions of $V$ such that
$\frac{\alpha_{V_i }}{\alpha_G} \ge c$ \,  $i\in [q]$, or equivalently, 
$ c\le \frac{\alpha_{V_i }}{\alpha_{V_j }} \le \frac1{c}$ $(i\ne j)$.

A more accurate balancing is defined by fixing a vector
$\a =( a_1 , \dots ,a_q )$ with components forming  a probability 
distribution over $[q]$: $a_i >0$ $i\in [q]$, $\sum_{i=1}^q a_i =1$. 
Let ${\cal P}^{\a}_q$ denote the set of $q$-partitions of $V$ such that 
$
 \left( \frac{\alpha_{V_1}}{\alpha_G} ,\dots ,\frac{\alpha_{V_q}}{\alpha_G}
 \right)
$ 
is approximately $\a$-distributed, that is
$$
 \left| \frac{\alpha_{V_i}}{\alpha_G} -a_i \right| \le 
   \frac{\alpha_{\max} (G)}{\alpha_G} \quad (i=1,\dots ,q),
$$
the right hand side tending to 0 as $|V(G)|\to\infty$ for 
weighted graphs with  no dominant vertex-weights. 

The {\it microcanonical ground state energy} of $G$ given $\a$ and $\J$ 
($\h =\0$) is
$$
{\hat {\cal E}}_q^{\a} (G ,\J ) =-\max_{P\in {\cal P}_q^{\a} } 
 \sum_{i,j\in[q]} \alpha_i (G/P) \alpha_j (G/P) \beta_{ij} (G/P) J_{ij} .
$$

\begin{rem}\label{energia1}
In Theorem 2.14 of~\cite{LovII} it is proved that the convergence of the
weighted graph sequence $(G_n )$ with no dominant vertex-weights is equivalent
to the convergence of its microcanonical ground state energies for any $q$, 
$\a$,
and $\J$. Also, it is equivalent to the convergence of its $q$-quotients in
Hausdorff distance for any $q$.
\end{rem}

\begin{rem}\label{energia2}
 Under the same conditions, Theorem 2.15 of~\cite{LovII}
states that the convergence of the above $(G_n )$ implies the convergence
of its ground state energies for any $q$, $\J$, and $\h$; further the
convergence of the spectrum of $(G_n )$.  
\end{rem}

Using these facts, we investigate the testability of some special 
multiway cut densities defined in the forthcoming definitions.

\begin{defn}
The {\it minimum q-way cut density} of $G$ is
$$
f_q (G) =
 \min_{P \in {\cal P}_q} \,
 \frac1{\alpha_{G}^2} \sum_{i=1}^{q-1} \, \sum_{j=i+1}^{q} e_G (V_i, V_j )  ,
$$ 
the {\it minimum  c-balanced q-way cut density} of $G$ is
\begin{equation}\label{fc}
f^c_q (G) =
 \min_{P \in {\cal P}^c_q} \,
 \frac1{\alpha_G^2} \, \sum_{i=1}^{q-1} \sum_{j=i+1}^{q} e_G (V_i, V_j )  ,
\end{equation}
and the {\it minimum}  $\a$-{\it balanced q-way cut density} of $G$ is
$$ 
f^{\a}_q (G) =
 \min_{P \in {\cal P}^{\a}_q} \,
 \frac1{\alpha_G^2} \, \sum_{i=1}^{q-1} \sum_{j=i+1}^{q} e_G (V_i, V_j )  .
$$ 
\end{defn}

Occasionally, we want to penalize cluster volumes that wildly differ.
For this purpose we herein introduce the notions of weighted 
minimum cut densities.

\begin{defn}
The {\it minimum weighted  q-way cut density} of $G$ is
$$
 \mu_q (G) = \min_{P \in {\cal P}_q} \,
 \sum_{i=1}^{q-1} \sum_{j=i+1}^{q}  \frac1{\alpha_{V_i } \cdot \alpha_{V_j } } 
  \cdot e_G (V_i, V_j ) , 
$$
the {\it minimum weighted  c-balanced q-way cut density} of $G$ is
$$
 \mu^c_q (G) = \min_{P \in {\cal P}^c_q} \,
 \sum_{i=1}^{q-1} \sum_{j=i+1}^{q}  \frac1{\alpha_{V_i }\cdot \alpha_{V_j } }
 \cdot e_G (V_i, V_j ) ,
$$
and the {\it minimum weighted}
 $\a$-{\it balanced q-way cut density}  of $G$ is
$$
 \mu^{\a}_q (G) = \min_{P \in {\cal P}^c_q} \,
 \sum_{i=1}^{q-1} \sum_{j=i+1}^{q}  \frac1{\alpha_{V_i }\cdot \alpha_{V_j } }
 \cdot e_G (V_i, V_j ) .
$$
\end{defn}

\begin{prop}\label{p1}
$f_q (G)$ is testable for any $q\le |V(G)|$.
\end{prop}

\begin{pf}
Observe that $f_q (G)$ is a special ground state energy:
$$
f_q (G) = {\hat {\cal E}}_{q} (G,\J , \0 ) ,
$$ 
where the magnetic field is $\0$ and the $q\times q$ symmetric matrix
$\J$ is the following:
$J_{ii}=0$ $i\in [q]$, further $J_{ij}=-1/2$ $(i\ne j)$.
By Remark~\ref{energia1}  and the  equivalent statement (c) of 
Theorem~\ref{ekvi},
the minimum $q$-way cut density is testable for any $q$.
\end{pf}

However, this statement is of not much use, since $f_q (G_n) \to 0$,
in the lack of dominant vertex-weights. In fact, the minimum $q$-way
cut density is trivially estimated from above by
$$
 f_q (G_n) \le (q-1) \frac{\alpha_{max} (G_n )}{\alpha_{G_n }} 
 + {{q-1} \choose 2} \left( \frac{\alpha_{max} (G_n )}{\alpha_{G_n }} \right)^2
$$
that tends to 0 provided ${\alpha_{\max} (G_n )}/{\alpha_{G_n }} \to 0$ as
$n\to\infty$.

\begin{prop}
$f_q^{\a} (G)$ is testable for any  $q\le |V(G)|$ and distribution $\a$
over $[q]$.
\end{prop}

\begin{pf}
Choose $\J$ as in the proof of Proposition~\ref{p1}. In this way, 
$f_q^{\a} (G)$ is a special microcanonical ground state energy:
\begin{equation}\label{fa}
f_q^{\a} (G) = {\hat {\cal E}}_q^{\a} (G,\J ) .
\end{equation}
Hence, 
by Remark~\ref{energia1}, the convergence of $(G_n)$ is equivalent to the
convergence of $f_q^{\a} (G_n)$ for any $q$ and 
 any distribution $\a$ over $[q]$. Therefore, by 
the equivalent statement (c) of Theorem~\ref{ekvi}, the testability of the
minimum $\a$-balanced $q$-way cut density also follows.
\end{pf}

\begin{prop}\label{fuc}
$f_q^{c} (G)$ is testable for any  $q\le |V(G)|$ and $c\le 1/q$.
\end{prop}

\begin{pf}
Theorem 4.7 and Theorem 5.5 of~\cite{LovII} imply that 
for  any two weighted graphs $G,G'$ 
\begin{equation}\label{ea}
  |{\hat {\cal E}}_{q}^{\a} (G,\J ) -{\hat {\cal E}}_{q}^{\a} (G',\J ) | 
\le (3/2 +\kappa ) \cdot d^{\text{Hf}}({\hat {S}}_q (G),{\hat {S}}_q (G') ) ,
\end{equation}
where $\kappa =o (\min \{ |V(G)|,|V(G')| \} )$ is a negligible small constant,
provided  the number of vertices of $G$ and $G'$ is sufficiently large.
By Remark~\ref{energia2} 
we know that if $(G_n)$ converges, its
$q$-quotients also converge in Hausdorff distance, consequently form a 
Cauchy-sequence. This means that for any $\eps >0$ there is an $N_0$
such that for $n,m>N_0$:  $d^{\text{Hf}} ( {\hat {S}}_q (G_n), 
{\hat {S}}_q (G_m) ) <\eps$.
We want to prove that for $n,m>N_0$:
$|f_q^c (G_n ) - f_q^c (G_m )| < 2\eps$. 
On the contrary, suppose that there are $n,m>N_0$ such that
$|f_q^c (G_n ) - f_q^c (G_m )|\ge 2\eps$. 
Say,
$f_q^c (G_n ) \ge f_q^c (G_m )$.  
Let $A:=\{ \a \, : a_i \ge c \, , \, i=1,\dots ,q \}$ is the subset of special 
$c$-balanced distributions  over $[q]$. On the one hand, 
$$f_q^c (G_m ) = \min_{\a \in A} f_q^{\a } (G_m ) =f_q^{\a^*} (G_m )$$
for some $\a^* \in A$. 
On the other hand, by~(\ref{fa}) and (\ref{ea}), 
$f_q^{\a^* } (G_n ) - f_q^{\a^* } (G_m ) \le (\frac32 +\kappa ) \eps$,
that together with the indirect assumption implies that
 $f_q^c (G_n ) - f_q^{\a^* } (G_n ) \ge (\frac12 -\kappa ) \eps >0$ for this
$\a^* \in A$. But this contradicts to the fact that $f_q^c (G_n )$ is the
minimum of $f_q^{\a } (G_n )$'s over $A$.
Thus,  $f_q^c (G_n )$ is also a Cauchy sequence, and being a real sequence,
it is also convergent.
\end{pf}

Concerning  the penalized densities, trivially,
$$
 \mu_q (G) 
 =\min_{P \in {\cal P}_q}  \sum_{i=1}^{q-1} \sum_{j=i+1}^{q} 
 \beta_{ij} (G/P) .
$$
In fact, $\mu_q (G)$ is not testable as
we can show  an example where $\mu_q (G_n )\to 0$,
but randomizing a sufficiently large part of $G_n$, the weighted
minimum $q$-way cut density of that part  is constant.
The example is for $q=2$ and for a simple graph on $n$ vertices
such that order of
$\sqrt{n}$ vertices are connected with a single edge to the remaining
vertices that form a complete graph. Then $\mu_2 (G_n)\to 0$, but
randomizing a sufficiently large part of the graph, with high probability,
it will be a subgraph of the complete graph, whose minimum 2-way cut density
is of constant order.

\begin{prop}\label{p2}
$\mu_q^{\a} (G)$ is testable for any  $q\le |V(G)|$ and distribution $\a$
over $[q]$.
\end{prop}

\begin{pf}
By  the definition of Hausdorff distance, the convergence of
 $q$-quotients  guarantees  the convergence of
\begin{equation}\label{mu}
 \mu_q^{\a} (G) =
 \min_{P\in {\cal P}_q^{\a}}  \sum_{i=1}^{q-1} \sum_{j=i+1}^{q} 
 \beta_{ij} (G/P) 
\end{equation}
for any $\a$ and $q$ in the following way. 
Let ${\hat {\cal S}}_q^{\a} (G)$ denote the set of factor graphs of $G$ 
with respect to
partitions in  ${\cal P}_q^{\a}$. As a consequence of
Lemma 4.5 and Theorem 5.4 of~\cite{LovII},
for any two weighted graphs $G$, $G'$ 
\begin{equation}\label{becsles}
 \max_{\a} d^{Hf} ({\hat {\cal S}}_q^{\a } (G) , 
                   {\hat {\cal S}}_q^{\a } (G') ) \le (3 +\kappa )\cdot 
  d^{Hf} ({\hat {\cal S}}_q (G) , {\hat {\cal S}}_q (G') ) ,
\end{equation}
where $\kappa =o (\min \{ |V(G)|,|V(G')| \} )$.

By Remark~\ref{energia1},
for a convergent graph-sequence $(G_n)$, the sequence 
${\hat {\cal S}}_q (G_n)$ converges,  and by the inequality~(\ref{becsles}), 
${\hat {\cal S}}_q^{\a } (G_n )$ also converges
in  Hausdorff distance for any distribution $\a$ over $[q]$.
As they form a Cauchy sequence,
 $\forall \eps$ $\exists N_0$ 
such that for $n,m>N_0$ 
$$
 d^{Hf} ({\hat {\cal S}}_q^{\a } (G_n ) , 
         {\hat {\cal S}}_q^{\a } (G_m ) ) <\eps 
$$
uniformly for any $\a$.
In view of the Hausdorff distance's definition, this means that for any
$q$-quotient $H\in {\hat {\cal S}}_q^{\a } (G_n )$ there exists (at least one) 
$q$-quotient $H' \in {\hat {\cal S}}_q^{\a } (G_m )$, and vice versa, for any
$H' \in {\hat {\cal S}}_q^{\a } (G_m )$ there exists (at least one)
$H \in {\hat {\cal S}}_q^{\a } (G_n )$ such that  $d_1 (H,H' ) <\eps$. 
(In fact, the maximum distance between the elements of the above pairs
is less than $\eps$. Note that the symmetry in the definition of the
Hausdorff distance is important: the pairing exhausts the sets
even if they have different cardinalities.)

Using the fact that the vertex-weights of such a pair
$H$ es $H'$ are almost the same (the coordinates of the vector $\a$),
by the notation $a=\min_{i\in [q]} a_i$, the following argument is valid
for $n,m$ large enough:
\begin{equation}
\begin{aligned}
&2a^2 \sum_{i\ne j} |\beta_{ij} (H) -\beta_{ij} (H')| \le
   \sum_{i,j=1}^q a^2 |\beta_{ij} (H) -\beta_{ij} (H')| \le \\
 &\le \sum_{i,j=1}^q | \alpha_i (H)
  {\alpha_j (H)} \beta_{ij} (H) -
  \alpha_i (H') \alpha_j (H')\beta_{ij} (H')| =d_1 (H,H')<\eps .
\end{aligned}
\end{equation}
Therefore
$$
 | \sum_{i=1}^{q-1} \sum_{j=i+1}^{q} \beta_{ij} (H) -
   \sum_{i=1}^{q-1} \sum_{j=i+1}^{q} \beta_{ij} (H') | < \frac{\eps}{2a^2}
  := {\eps}' ,
$$
and because 
$\sum_{i=1}^{q-1} \sum_{j=i+1}^{q} \beta_{ij} (H)$ and
$\sum_{i=1}^{q-1} \sum_{j=i+1}^{q} \beta_{ij} (H')$
are individual terms behind the minimum in~(\ref{mu}), the above inequality
holds for their minima over ${\cal P}_q^{\a}$ as well:
\begin{equation}\label{szam}
 | \mu_q^{\a} (G_n) - \mu_q^{\a} (G_m) | <  {\eps}' .
\end{equation}
Consequently, the sequence $\mu_q^{\a} (G_n)$ is a Cauchy sequence, and
being a real sequence, it is also convergent. Thus $\mu_q^{\a}$ is testable.
\end{pf}

\begin{rem}
The testability of $\mu_q$, apparently, does not follow in the same way
due to presence of distinct  vetex-weights in $H$ and $H'$. Thus, the smallness
of $d_1 (H,H')$ does not imply the closeness of their edge-weights.
\end{rem}

However, as the testability of $f_q^{\a}$  implied  the
testability of $f_q^{c}$,   the testability of
$\mu_q^{\a}$  also implies the  testability  of
$\mu_q^c$.

\begin{prop}\label{muc}
$\mu_q^{c} (G)$ is testable for any  $q\le |V(G)|$ and $c\le 1/q$.
\end{prop}

\begin{pf}
The proof is analogous to that of Proposition~\ref{fuc} using
equation~(\ref{szam}) instead of equation~(\ref{ea}).
By the pairing argument of the proof of Proposition~\ref{p2}, the real sequence
$\mu^c_q (G_n )$ is a Cauchy sequence, and therefore, convergent.
This immediately implies the testability of $\mu^c_q $.
\end{pf}

By Remark~\ref{energia2},
 the convergence of $(G_n )$ also implies the convergence of the
spectra, though the convergence of the spectrum itself is weaker than the 
convergence of the graph sequence. Without going into details, we remark that
in~\cite{Bolla}, $f_q (G)$ and $\mu_q (G)$ were bounded from below by the 
$q$ smallest Laplacian eigenvalues of $G$. An upper estimate can also
be constructed and we conjecture that in case of testable parameters an
asymptotic estimate is also valid. 


\section{Minimum cut as a quadratic programming problem}\label{prob}

In Section 4, we proved that $f^c_q$  is a testable weighted  graph parameter. 
Now, we extend it to graphons.

\begin{prop}\label{ext}
Let us define the graphon functional ${\tilde{f}}^c_q$ in the following way:
\begin{equation}\label{inf}
 {\tilde{f}}^c_q (W ) := \inf_{Q \in {\cal Q}^c_q } 
   \sum_{i=1}^{q-1} \sum_{j=i+1}^{q} 
 \iint_{ S_i \times S_j } W(x,y) \, dx \, dy   =
 \inf_{Q \in {\cal Q}^c_q }   {\tilde{f}}_q (W; S_1 ,\dots ,S_q )
\end{equation}
where the infimum is taken over all the $c$-balanced Lebesgue-measurable
partitions
$Q =(S_1 ,\dots ,S_q)$ of
[0,1]. For these, $\sum_{i=1}^q \lambda (S_i ) =1$ and
$\lambda (S_i ) \ge c$ \, ($i\in [q]$), where $\lambda$ denotes the
Lebesgue-measure, and
${\cal Q}^c_q $ denotes the set of $c$-balanced $q$-partitions of [0,1].
We state that ${\tilde{f}}^c_q$ is the extension of ${f}^c_q$ in the
following sense: If $(G_n)$ is a convergent weighted  graph sequence with
uniformly bounded edge-weights and no dominant vertex-weights, then
denoting by $W$ the
essentially unique limit graphon of the sequence (see Theorem~\ref{lim}), 
${f}^c_q (G_n) \to {\tilde{f}}^c_q (W)$ as $n\to\infty$.
\end{prop}

\begin{pf}
First we show that  ${\tilde{f}}^c_q $ is continuous in the cut-norm.
As ${\tilde{f}}^c_q (W)$ is insensitive to measure preserving maps of $W$, 
it suffices to prove that
to any $\eps$ we can find $\eps'$ such that for any two
graphons $W,U$ with $\| W -U \|_{\square} <\eps'$, the relation
$|{\tilde{f}}^c_q (W )- {\tilde{f}}^c_q (U )|<\eps $ also holds.
By the definition of the cut-norm, for any Lebesgue-measurable
$q$-partition $(S_1 ,\dots ,S_q )$ of [0,1], the relation
$$
 |\iint_{ S_i \times S_j } (W(x,y)   -
 U (x,y) ) \, dx \, dy  | \le \eps'  \qquad (i\ne j) 
$$
holds. Summing up for the $i\ne j$ pairs
\begin{equation}\label{est}
 | \sum_{i=1}^{q-1} \sum_{j=i+1}^q \iint_{ S_i \times S_j } W(x,y) \, dx \, dy
-\sum_{i=1}^{q-1} \sum_{j=i+1}^q \iint_{ S_i \times S_j } U (x,y) 
\, dx \, dy | \le {q \choose 2} \eps' .
\end{equation}
Therefore
$$
\inf_{(S_1 ,\dots ,S_q) \in {\cal Q}^c_q } 
   \sum_{i=1}^{q-1} \sum_{j=i+1}^{q} 
 \iint_{ S_i \times S_j } W(x,y) \, dx \, dy  \ge
\inf_{(S_1 ,\dots ,S_q) \in {\cal Q}^c_q } 
   \sum_{i=1}^{q-1} \sum_{j=i+1}^{q} 
 \iint_{ S_i \times S_j } U(x,y) \, dx \, dy - {q \choose 2} \eps' 
$$
and vice versa,
$$
\inf_{(S_1 ,\dots ,S_q) \in {\cal Q}^c_q } 
   \sum_{i=1}^{q-1} \sum_{j=i+1}^{q} 
 \iint_{ S_i \times S_j } U(x,y) \, dx \, dy  \ge
\inf_{(S_1 ,\dots ,S_q) \in {\cal Q}^c_q } 
   \sum_{i=1}^{q-1} \sum_{j=i+1}^{q} 
 \iint_{ S_i \times S_j } W(x,y) \, dx \, dy - {q \choose 2} \eps'  .
$$
Consequently the absolute difference of the two infima is bounded from
above by  ${q \choose 2} \eps'$.
Thus, $\eps' =\eps /{q \choose 2}$ will do.


Let $(G_n)$ be a convergent weighted  graph sequence with
uniforly bounded edge-weights and no dominant vertex-weights. By
Theorem~\ref{lim},
there is an essentially unique graphon $W$ such that $G_n \to W$,  i.e., 
$\delta_{\square} (W_{G_n} , W) \to 0$ as $n\to\infty$. By the continuity of
 ${\tilde{f}}^c_q$, 
\begin{equation}\label{tart}
{\tilde{f}}^c_q (W_{G_n }) \to {\tilde{f}}^c_q (W), \qquad n\to\infty .
\end{equation} 
Suppose that
$$
 {\tilde{f}}^c_q (W_{G_n} ) ={\tilde{f}}_q (W_{G_n}; S_1^* ,\dots ,S_q^* ) ,
$$
that is the infimum in~(\ref{inf}) is attained at the 
$c$-balanced Lebesgue-measurable
$q$-partition $(S_1^* ,\dots ,S_q^* )$ of [0,1].

Let $G_{nq}^*$ be the $q$-fold blown-up  of $G_n$ with respect to
$(S_1^* ,\dots ,S_q^* )$. It is a weighted graph on at most $nq$
vertices  defined in the following way. 
Let $I_1 ,\dots ,I_n$ be  consecutive intervals of [0,1] such that
$\lambda (I_j ) =\alpha_j (G_n )$, $j=1,\dots ,n$. 
The weight of the vertex labeled by $ju$ of
$G_{nq}^*$ is $\lambda (I_j \cap S_u^*)$,
$u\in[q]$, $j\in [n]$, while the edge-weights are 
$\beta_{ju,iv} (G_{nq}^*) = \beta_{ji} (G_n )$.
Trivially, the graphons $W_{G_n}$ and  $W_{G_{nq}^*}$ essentially define
the same stepfunction, hence 
${\tilde{f}}^c_q  (W_{G_n}) = {\tilde{f}}^c_q (W_{G_{nq}^*})$.
Therefore, by~(\ref{tart}), 
\begin{equation}\label{tart1}
 {\tilde f}^c_q (W_{G_{nq}^*}) \to {\tilde{f}}^c_q (W), \qquad n\to\infty .
\end{equation} 
As  $\delta_{\square} (G_n , G_{nq}^*) =
 \delta_{\square} (W_{G_n} , W_{G_{nq}^*})=0$, by part (e) of
Theorem~\ref{ekvi}  it follows that 
\begin{equation}\label{tart2}
 |{f}^c_q (G_{nq}^*) -  {f}^c_q (G_{n}) | \to 0, \qquad n\to\infty .
\end{equation}
Finally, by the construction of  $G_{nq}^*$,
${\tilde{f}}^c_q (W_{G_{nq}^*}) ={f}^c_q (G_{nq}^*)$, and hence,
$$
 |{f}^c_q (G_n) - {\tilde{f}}^c_q (W) | \le
 |{f}^c_q (G_n) - {f}^c_q (G_{nq}^*) | + 
 |{\tilde f}^c_q (W_{G_{nq}^*})- {\tilde{f}}^c_q (W)|
$$
that, in view of (\ref{tart1}), (\ref{tart2}), implies the required 
statement.
\end{pf}

\begin{cor}\label{kov}
In Section 3, while proving Theorem~\ref{ekvi}, 
an essentially unique extension of
a testable graph parameter to graphons was given. By  Proposition~\ref{ext},
the above ${\tilde{f}}^c_q$ is the desired extension 
of ${f}^c_q$, therefore
part (d) of Theorem~\ref{ekvi} is also applicable to it:
For a weighted graph sequence $(G_n )$ with 
$ \max_i \frac{\alpha_i (G_n ) }{\alpha_{G_n} } \to 0 $, the limit relation
$\tilde{f^c_q} (W_{G_n} ) -f^c_q(G_n) \to 0$ also holds as $n\to\infty$.
\end{cor}


Corollary~\ref{kov} gives rise to approximate the minimum $c$-balanced $q$-way
cut density  of a weighted graph on ``many'' vertices with no dominant
vertex weights by the extended $c$-balanced $q$-way cut density of the
 stepfunction graphon assigned to the graph.
In this way, the discrete optimization problem can be formulated as a
quadratic programming task with linear equality and inequality constraints.

To this end,
let us investigate a fixed weighted graph $G$ on $n$ vertices ($n$ is large).
To simplify notation we drop the subscript $n$, and $G$ in the arguments
of the vertex- and edge-weights. 
As $f_q^c (G)$ is invariant under the scale of the vertices, we can suppose
that $\alpha_G =\sum_{i=1}^n \alpha_i =1$. 
As $\beta_{ij}\in [0,1]$, $W_G$ is uniformly bounded by 1. Recall  that
$W_G (x,y) =\beta_{ij}$, if 
$x\in I_i$, $y \in I_j$, where $\lambda (I_j ) =\alpha_j$ $(j=1,\dots ,n)$
and  $I_1 ,\dots ,I_n$ are consecutive intervals of [0,1].

For fixed $q$ and $c\le 1/q$,
$f_q (G; V_1 ,\dots ,V_q) =
\frac1{\alpha_G^2} \, \sum_{i=1}^{q-1} \sum_{j=i+1}^{q} e_G (V_i, V_j )$
is a function taking on discrete values
over $c$-balanced $q$-partitions $P=(V_1 ,\dots ,V_q) \in
{\cal P}^c_q$ of the vertices of $G$.
As $n\to\infty$, by Corollary~\ref{kov}, this function approaches 
${\tilde{f}}_q (W_G ; S_1 ,\dots ,S_q )$ that is already a continuous 
function over $c$-balanced $q$-partitions
$Q =(S_1 ,\dots ,S_q) \in {\cal Q}_q^c$ of [0,1].
In fact, 
this continuous function can be regarded as a
multilinear function of the variable 
$$
 \x =(x_{11},\dots,x_{1n},x_{21},\dots, x_{2n},\dots ,x_{q1},\dots ,x_{qn})^T
 \in\R^{nq}
$$
where the coordinate indexed by $ij$ is
$$
 x_{ij} = \lambda (S_i \cap I_j ) , \quad j=1,\dots ,n; \quad i=1,\dots q.
$$
Hence, 
$$
 {\tilde f}_q (W_G; S_1 ,\dots ,S_q )=
 {\tilde f}_q (\x ) =\sum_{i=1}^{q-1} \sum_{i'=i+1}^q 
                     \sum_{j=1}^n \sum_{j'=1}^n x_{ij} x_{i' j'} \beta_{jj'}
  = \frac12 \x^T (\A \otimes \BB ) \x  ,
$$
where -- denoting by $\1_{q\times q}$ and $\I_{q\times q}$ the $q\times q$ 
all 1's and the identity matrix, respectively -- 
the eigenvalues of the $q\times q$ symmetric matrix $\A =\1_{q\times q} -
\I_{q\times q}$ are the number $q-1$ and -1 with multiplicity $q-1$, 
while those of the
$n\times n$ symmetric matrix $\BB =(\beta_{ij})$ are $\lambda_1 \ge \dots \ge
\lambda_n$. Latter one being a Frobenius-type matrix, $\lambda_1 >0$.
The eigenvalues of the Kronecker-product $\A \otimes \BB$ are the  numbers
$(q-1)\lambda_i$ $(i=1,\dots ,n )$ and $-\lambda_i$ with multiplicity $q-1$
$(i=1,\dots ,n)$. Therefore the above quadratic form is indefinite.
 
Hence, we have the following {\it quadratic programming} task:
\begin{equation}\label{dom}
\begin{aligned}
&\text{minimize} \qquad
{\tilde f}_q (\x ) = \frac12 \x^T (\A \otimes \BB ) \x  \\
&\text{subject to} \quad
\x \ge 0; \quad
\sum_{i=1}^q x_{ij} = \alpha_j \quad (j\in[n]); \quad
 \sum_{j=1}^n x_{ij}  \ge c \quad (i\in[q]).
\end{aligned}
\end{equation}
The feasible region is the closed convex polytope of~(\ref{dom}), and it is,
in fact, in an $n(q-1)$-dimensional hyperplane of $\R^{nq}$. The gradient
of the objective function
$\nabla {\tilde f}_q (\x ) =(\A \otimes \BB ) \x $ cannot be $\0$
in the feasible region, 
 provided the weight matrix $\BB$, and hence $\A \otimes \BB$ is
non singular.

The arg-min of the quadratic programming task~(\ref{dom}) is one of the
Kuhn--Tucker points (giving relative minima of the indefinite quadratic
form over the feasible region), that can be found by numerical algorithms
(by tracing back the problem to a linear programming task), see~\cite{Baz}.

Eventually,
we give the extension of the testable weighted graph parameter $\mu_q^c$
to graphons.

\begin{prop}
Let us define the graphon functional ${\tilde{\mu}}^c_q$ in the following way:
$$
 {\tilde{\mu }}^c_q (W ) := \inf_{Q \in {\cal Q}^c_q } 
   \sum_{i=1}^{q-1} \sum_{j=1}^{q} 
 \iint_{ S_i \times S_j } \frac1{\lambda (S_i) \lambda(S_j)}\, 
 W(x,y)\, dx\, dy .
$$
We state that ${\tilde{\mu }}^c_q$ is the extension of ${\mu }^c_q$ in the
following sense: If $(G_n)$ is a convergent weighted  graph sequence with
uniforly bounded edge-weights and no dominant vertex-weights, then
denoting by $W$ the
essentially unique limit graphon of the sequence (see Theorem~\ref{lim}), 
${\mu }^c_q (G_n) \to {\tilde{\mu }}^c_q (W)$ as $n\to\infty$.
\end{prop}

The proof is analogous to that  of Proposition~\ref{ext}, after we 
have proved that
 ${\tilde{\mu }}^c_q $ is continuous in the cut-norm.
In fact, with estimates, analogous to~(\ref{est}),
 $\eps' =\eps c^2/{q \choose 2}$ will do.

Consequently,
${\tilde{\mu }}^c_q (W_G ) -\mu_q^c (G) \to 0$
as $V(G)\to \infty$ with no dominant vertex-weights.
This fact also gives rise to approximate the minimum $c$-balanced 
weighted $q$-way cut density of a large graph by quadratic programming methods.

\section{Convergence of noisy graph sequences}\label{noise}

Now, we use the above theory for perturbations.
If not stated otherwise, the vertex-weights are equal (say 1), 
and a weighted graph $G$ on $n$
vertices is identified with its $n\times n$ symmetric weight matrix $\A$.
Let $G_{\A}$ denote the weighted graph with unit vertex-weights and
edge-weights that are entries of  $\A$ .

\begin{defn}
Let $w_{ij}$ $(1\le i \le j\le n)$
be independent random variables  defined on the same
probability space, and $w_{ji}=w_{ij}$. $\EE (w_{ij}) =0$  $(\forall i,j)$
and the $w_{ij}$'s  are uniformly bounded,
i.e., there is a constant $K>0$ -- that does not depend of $n$ -- 
 such that $|w_{ij}| \le K$, $\forall i,j$.
The  $n\times n$ symmetric real random matrix 
$\W =(w_{ij})_{1\le i\le n, \, 1\le j\le n }$ 
is called a Wigner-noise.
\end{defn}

The corresponding edge-weighted graph $G_{\W}$ is called a {\it Wigner-graph}.
To indicate that the size $n$ is expanding, we use the notations
$\W_n$ and $G_{\W_n}$.

\begin{defn}
The $n\times n$ symmetric real matrix $\BB$ is a 
blown-up  matrix,
if there is a 
$q\times q$ symmetric so-called pattern matrix
$\PPP$ with entries $0 < p_{ij} < 1$, 
and there  are positive integers 
$n_1 ,\dots ,n_q$ with $\sum_{i=1}^q n_i =n$,
such that -- after rearranging its rows and columns --
the matrix $\BB$ can be divided into $q\times q$ blocks, where 
block $(i,j)$ is 
an $n_i \times n_j$ matrix with entries 
all equal to $p_{ij}$ $(1\le i,j \le n )$.
\end{defn}

Fix $\PPP$,  blow it up to an $n\times n$ matrix  $\BB_n$, and consider
the noisy matrix $\A_n =\BB_n +\W_n$ as $n_1 , \dots ,n_q \to\infty$
at the same rate. 

\begin{rem}\label{K}
While perturbing $\BB_n$ by $\W_n$, for the uniform bound of the entries of 
$\W_n$ the condition
\begin{equation}\label{KK}
K\le \min \{ \min_{i,j\in [q]} p_{ij}\, , \, 1-\max_{i,j\in [q]} p_{ij} \}
\end{equation}
is satisfied. In this way, the entries of $\A_n$  are in the [0,1]
interval, and hence, $G_{\A_n} \in {\cal G}$. 

We remark that $G_{\W_n}
\notin {\cal G}$, but $W_{G_{\W_n}} \in {\cal W}$ and the theory of bounded
graphons applies to it. By adding an appropriate Wigner-noise to $\BB_n$,
we can achieve that  $\A_n$ becomes a 0-1 matrix:
its entries are equal to 1 with probability $p_{ij}$ and 0 otherwise
within the block of size $n_i \times n_j$ (after rearranging its rows and
columns). In this case, the corresponding noisy graph $G_{\A_n}$ is
a random simple graph.
\end{rem}

As $\rk (\BB_n ) =q$ and $\| \W_n \| ={\cal O} (\sqrt{n} )$ almost surely
($n\to\infty$), 
 the noisy matrix $\A_n$ almost surely has $q$ protruding eigenvalues 
(of order $n$), and all the other eigenvalues are of order $\sqrt{n}$,
there is a  spectral gap between the $q$ largest and the
other eigenvalues $\A_n$. 

Let $\X_n = (\x_1 ,\dots ,\x_q)$ be the
  $n\times q$ matrix containing the eigenvectors belonging to the $q$
protruding eigenvalues of $\A_n$ in its columns.
The rows of  $\X_n$, that is the vectors $\x^{1}, \dots ,\x^{n} \in \R^q $ 
are regarded as $q$-dimensional representatives of the vertices of $G_{\A_n}$.
The $q$-variance of the representatives is
$$
 S_q^2 (\X_n )= 
 \sum_{i=1}^q \sum_{j\in V_i} \|\x^{j}- {\bar \x}^{i}\|^2 ,
$$
where ${\bar \x}^{i} =\frac1{n_i} \sum_{j\in V_i} \x^{j}$. 

In~\cite{Bolla2} we proved that
$$  
  S_q^2 (\X_n ) = {\cal O} \left( \frac1{n} \right)
$$
almost surely, under the growth condition $n_i /n \ge c$ $(i=1,\dots ,q)$.

In the other direction: 
for sufficiently large $n$, under some conditions, we can
separate an $n\times n$ symmetric ``error-matrix'' $\EEE$ from $\A$,
such that $\| \EEE \| ={\cal O} (\sqrt{n})$ and the remaining matrix
$\A -\EEE$  is a blown-up matrix $\BB$ of ``low rank''.
Consequently, 
$G_{\BB }$ is a weighted graph with homogeneous edge-densities within
the clusters (determined by the blow-up).
It resembles to the weak Szemer\'edi-partition, cf.~\cite{LovSzeg2}, 
but the error-term is bounded in spectral norm, instead of the cut-norm.
However, by large deviations, we can prove that  the cut-norm of a 
Wigner-graph tends to zero
almost surely as  $n\to \infty$.

\begin{thm}\label{wign0}
 For  any sequence $(G_{\W_{n}})$  of  Wigner-graphs
$$
\lim_{n\to \infty} \| W_{G_{\W_n}} \|_{\square} =0 \qquad (n\to\infty )
$$
almost surely.
\end{thm}

To prove the theorem, we need a proposition
that is an easy consequence of Azuma's martingale inequality,
see Theorem 5.3 of~\cite{Chung}.
 
\begin{prop}\label{azuma}
Let
 $X_1, X_2, \dots, X_N$ be i.i.d. random variables
 with zero mean and $|X_i |\le 1$, $i=1,\dots ,n$. Then
\begin{equation}
\PP\left( \left \vert \sum_{j=1}^N X_j\right \vert > \gamma N \right)< 
2
 \exp\left (-\frac {N \gamma^2}{2}\right) , \qquad  0<\gamma <1 .
\label{ch1}
\end{equation}
\end{prop}

\begin{pf} Now we are ready to prove Theorem~\ref{wign0}.
By  the definition of the cut-norm of a stepfunction
graphon and \cite{LovSzeg1}, 
\begin{equation}
 \| W_{G_{\W_{n}}} \|_{\square}
 = \frac 1 {n^2} \max_{U, T\subset [n] }
 \left | \sum_{i\in U}  \sum_{j\in T} w_{ij} \right |.
\label{cutn}
\end{equation}
We remark that $\max_{U, T\subset [n] }
 \left | \sum_{i\in U}  \sum_{j\in T} w_{ij} \right |$ is the cut-norm
of the matrix $\W_n$ defined in~\cite{Frieze}.

To make the entries behind the double sum of (\ref{cutn}) independent, we use 
formulas (7.2), (7.3) of~\cite{LovI}:
\begin{equation}
 \| W_{G_{\W_n}} \|_{\square} \le 6 \max_{U\subset [n]}
 {1 \over n^2}\left |
 \sum_{i\in U}  \sum_{j\in [n] \setminus U }
w_{ij}
 \right | . 
 \label{cutfg} 
\end{equation}

Apply Proposition~\ref{azuma}
for a subsequence of length   $N$
 of entries of $\W_n$ which does not contain  
 $w_{ij}$ and $w_{ji}$ simultaneously.  
Namely,  $i\in U$, $j\in [n]\setminus U$, 
$N=|U|\cdot (n-|U| )$. Remark that $n-1\le N\le n^2 /4$.

We distinguish between two cases.
\begin{itemize}
\item {\it Case 1.}
Suppose that $N\le n^{3/2}$ and apply Proposition~\ref{azuma}
for the right 
 hand side of (\ref{cutfg}) with $\gamma = n^{4/10}$:
\begin{equation}\label{ch4}
\begin{aligned}
  & \PP\left( {1
 \over n^2}
\left |
 \sum_{i\in U}  \sum_{j\in [n] \setminus U }
w_{ij}
 \right |>
  n^{-1/10}  \right)
  =  \PP\left( 
\left |
 \sum_{i\in U}  \sum_{j\in [n] \setminus U }
w_{ij}
 \right |> n^{3/2}n^{4/10}
    \right)\le\\
 \le &\PP\left( 
\left |
 \sum_{i\in U}  \sum_{j\in [n] \setminus U }
w_{ij}
 \right |> N  n^{4/10}
 \right)
 < 2 \exp\left( - {n^{9/5} \over 2} \right).
\end{aligned}
\end{equation}
In the first inequality we used the condition  $N\le n^{3/2}$,
 while  the second inequality follows from (\ref{ch1}) and  the 
  fact that $N \ge n-1$. 

\item {\it Case 2.}
 Now suppose that $N > n^{3/2}$ and apply Proposition~\ref{azuma}
 for the right hand side of (\ref{cutfg}) with $\gamma = 4n^{-1/10}$:
\begin{equation}\label{ch5}
\begin{aligned}
  & \PP\left( {1
 \over n^2}
\left |
 \sum_{i\in U}  \sum_{j\in [n] \setminus U }
w_{ij}
 \right |>
  n^{-1/10}  \right)
  =  \PP\left( 
\left |
 \sum_{i\in U}  \sum_{j\in [n] \setminus U }
w_{ij}
 \right |> n^{2}n^{-1/10}
    \right)\le \\
 \le &\PP\left( 
\left |
 \sum_{i\in U}  \sum_{j\in [n] \setminus U }
w_{ij}
 \right |>4 N  n^{-1/10}
 \right)
 <2 \exp\left( - 8 n^{13/10} \right).
\end{aligned}
\end{equation}
In the first inequality we used the fact  $4N\le n^{2}$,
 while  the second inequality follows from (\ref{ch1}) and  the 
  condition  $N > n^{3/2}$. 
\end{itemize}

As for large values of $n$ the right hand side of (\ref{ch5}) is 
greater than that of (\ref{ch4}),
the  probability
$\PP\left( \max_{U\subset [n]} \left \vert \sum_{i\in U}\sum_{j \in 
[n]\setminus U} 
w_{ij}
\right \vert > n^{-1/10}\right)$
can be bounded by  the number of possible 2-partitions of [n] times
 the right hand side of (\ref{ch5}):
\begin{equation}
\PP\left( \| W_{G_{\W_n}} \|_{\square} 
>6\cdot (n)^{-1/10} \right)<
  2^{n+1}
 \exp\left(- 8 n^{13/10} 
 \right) .
 \label{chcut}
\end{equation}

 As the right hand side of  
(\ref{chcut}) 
 is a  general term of 
 a  convergent series,  the statement of the theorem follows
 by the  Borel-Cantelli Lemma. 
\end{pf}

\begin{rem}\label{zajostart}
Let $\A_n :=\BB_n +\W_n$ and $n_1 , \dots ,n_q \to\infty$ in such a way
that $\lim_{n\to\infty} \frac{n_i}{n} =r_i$ $(i=1,\dots ,q)$, 
$n=\sum_{i=1}^q n_i$; further, for the uniform bound $K$ of the entries of
the ``noise'' matrix $\W_n$ the condition~(\ref{KK}) is satisfied.
Under these conditions, Theorem~\ref{wign0} implies that the 
``noisy'' graph sequence  $(G_{\A_n })\subset {\cal G}$  
converges almost surely in
the $\delta_{\square}$ metric. It is easy to see that the almost sure limit is
the stepfunction $W_H$, where the factor graph $H=G_{\BB_n} /P$  
does not depend on $n$, as
$P$ is the $q$-partition of the vertices of  $G_{\BB_n}$ with resepect  to
the blow-up  (with cluster sizes $n_1 ,\dots ,n_q$). Actually,
the  vertex- and edge-weights of the weighted graph $H$ are
$$ 
  \alpha_i (H) =r_i \quad (i\in [q]), \qquad
  \beta_{ij} (H) = \frac{n_i n_j p_{ij}}{n_i n_j } =p_{ij} \quad (i,j\in [q]).
$$
\end{rem}

\begin{rem}
Under the conditions of Remark~\ref{zajostart},
as $(G_{\A_n }) \subset {\cal G}$ converges almost surely and 
$f_q$, $f_q^c$, $f_q^{\a}$,  $\mu_q^c$,  $\mu_q^{\a}$
are testable graph parameters,  the sequences
$f_q (G_{\A_n })$, $f_q^c (G_{\A_n })$, $f_q^{\a} (G_{\A_n })$, 
 $\mu_q^c (G_{\A_n })$,  $\mu_q^{\a} (G_{\A_n })$ also converge
almost surely. However, the almost sure limits are
the corresponding extended ${\tilde f}_q$- or ${\tilde \mu}_q$-values of the
graphon $W_{H}$ and not the $f_q$- or $\mu_q$-values of $H$. 
For example, in Section 4, we have shown that $f_q (G_{\A_n }) \to 0$
$(n\to\infty )$, but
$f_q (H) = \sum_{i=1}^{q-1} \sum_{j=i+1}^q p_{ij} \ne 0$.
\end{rem}


\begin{ack}
We are indebted to L\'aszl\'o Lov\'asz and Katalin Friedl for
inspiring discussions.
\end{ack}

\end{document}